\documentclass[a4paper, BCOR=0.0mm, DIV=calc]{scrartcl}
\usepackage[T1]{fontenc}
\usepackage[utf8]{inputenc}
\usepackage[ngerman]{babel}
\usepackage[osf,sc]{mathpazo}
\usepackage{textcomp}
\usepackage{microtype}
\usepackage{amsmath, amssymb, amsfonts}
\KOMAoptions{twoside=false, twocolumn=false, headinclude=false, footinclude=false, mpinclude=false, pagesize=auto}
\recalctypearea
\begin{document}
\boldmath
\title{Spekulation über die Integralformel \boldmath$\int \frac{x^n \mathrm{d}x}{\sqrt{aa-2bx+cxx}}$\unboldmath, wo zugleich herausragende Beobachtungen über Kettenbrüche entstehen\footnote{
Originaltitel: "`Speculationes super formula integrali $\int{(x^n\mathrm{d}x)/\sqrt{aa-2bx+cxx}}$, ubi simul egregiae observationes circa fractiones continuas occurrunt"', erstmals publiziert in "`\textit{Acta Academiae Scientarum Imperialis Petropolitinae} 1782, 1786, pp. 62-84"', Nachdruck in "`\textit{Opera Omnia}: Series 1, Volume 18, pp. 244 - 264"', Eneström-Nummer E606, übersetzt von: Alexander Aycock, Textsatz: Artur Diener,  im Rahmen des Projektes "`Eulerkreis Mainz"' }}
\unboldmath
\author{Leonhard Euler}
\maketitle
\paragraph{§1}
Wir wollen vom einfachsten Fall, in dem $n=0$ ist, beginnen und wollen das Integral der Formel $\frac{\mathrm{d}x}{\sqrt{aa-2bx-cxx}}$ suchen, die für $x=\frac{b+z}{c}$ übergeht in diese $\frac{\mathrm{d}z}{\sqrt{aacc-bbc+czz}}$, wo zwei Fälle unterscheiden werden sollten, je nachdem ob $c$ entweder eine positive Größe oder eine negative war. Es sei daher zuerst $c=+ff$ und unsere Formel wird $\frac{\mathrm{d}z}{f\sqrt{aaff-bb+zz}}$ sein, dessen Integral $\frac{1}{f}\log{\frac{z + \sqrt{aaff-bb+zz}}{c}}$ ist und daher wird unser Integral
\[
\frac{1}{\sqrt{c}}\log{\frac{cx-b+\sqrt{aac-2bcx+ccxx}}{c}}
\]
sein, was, also so genommen, dass sie für $x=0$ verschwindet, liefert
\[
\frac{1}{\sqrt{c}}\log{\frac{cx-b+\sqrt{c(aa-2bx+cxx)}}{-b+a\sqrt{c}}}
\]
Aber wenn $c$ eine negative Konstante war, z.\,B. $c=-gg$, wird die durch $z$ ausgedrückte Integralformel $\frac{dz}{g\sqrt{aagg+bb-zz}}$ dessen Integral $\frac{1}{g}\arcsin{\frac{z}{\sqrt{aagg+bb}}} + C$ ist, wodurch das Integral so genommen, dass es für $x=0$ verschwindet,
\[
=\frac{1}{g}\arcsin{cx-b}{\sqrt{aagg+bb}} + \frac{1}{g}\arcsin{\frac{b}{\sqrt{aagg+bb}}}
\]
sein wird.
\paragraph{§2}
Es bezeichne nun $\Pi$ den Wert der Integralformel $\int \frac{dx}{\sqrt{aa-2bx+cxx}}$ so genommen, dass sie für $x=0$ verschwindet, entweder ist $c$ eine positive Größe oder eine negative, und wenn $c=+ff$ ist, wird wie wir gesehen haben 
\[
\Pi = \frac{1}{f}\log{\frac{ffx-b+f\sqrt{aa-2bx+ffx}}{af+b}}
\]
in dem anderen Fall aber, in dem $c=-gg$ ist, wird 
\[
\Pi = -\frac{1}{g}\arcsin{\frac{ggx+b}{\sqrt{aagg+bb}}} + \frac{1}{g}\arcsin{\frac{b}{\sqrt{aagg+bb}}}
\]
sein, oder durch Zusammenziehen beider Arcusfunktionen haben wird
\[
\Pi = \frac{1}{g}\arcsin{\frac{bg\sqrt{aa-2bx-ggxx}-abg-ag^3x}{aagg+bb}}
\] 
Da wir ja bald zeigen werden, dass die Integration der allgemeinen Formel $\int\frac{x^ndx}{\sqrt{aa-2bx+xx}}$, immer zu Fall $n=0$ zurückgeführt werden kann, wenn $n$ eine positive ganze Zahl war, können alle Integrale durch $\Pi$ ausgedrückt werden.
\paragraph{§3}
Gleich nach der Integration wollen wir der Variablen Größe $x$ einen konstanten Wert zuteilen, durch den die irrationale Formel $\sqrt{aa-2bx+xx}$ zu Null gemacht wird, was geschieht, wenn für $x = \frac{b\pm\sqrt{bb-aac}}{c}$genommen wird, und daher in zwei Fällen. Wir setzen für beide Fälle, dass die Funktion $\Pi$ übergeht in $\Delta$, sodass im Fall $c=ff$
\[
\Delta = \frac{1}{f}\log{\frac{\sqrt{bb-aaff}}{af-b}} = \frac{1}{f}\log{\sqrt{\frac{b+af}{b-af}}}
\]
ist, für den anderen Fall aber, für $c=-gg$, dass 
\[
\Delta = \frac{1}{g}\arcsin{\frac{ag\sqrt{bb+aagg}}{aagg+bb}} = \frac{1}{g}\arcsin{\frac{ag}{\sqrt{bb+aagg}}}
\]
ist. Wir werden den Wert $\Delta$ in den folgenden Fällen, bei denen der Wurzelausdruck $\sqrt{aa-2bx+cxx}$ verschwindet, hauptsächlich betrachtet haben.
\paragraph{§4}
Um nun zum folgenden Fall voranzuschreiten, wollen wir die Formel $s = \sqrt{aa-2bx+cxx} - a$ betrachten, dass sie natürlich für $x=0$ verschwindet, und da ja 
\[
\mathrm{d}s = \frac{-b\mathrm{d}x+cx\mathrm{d}x}{\sqrt{aa-2bx+cxx}}
\]
ist, wird wiederum durch Integrieren
\[
c\int\frac{x\mathrm{d}x}{\sqrt{aa-2bx+cxx}} = b\int\frac{\mathrm{d}x}{\sqrt{aa-2bx+cxx}} + s
\]
sein, woher wir berechnen
\[
\int\frac{x \mathrm{d}x}{\sqrt{aa-2bx+cxx}} = \frac{b}{c}\Pi + \frac{\sqrt{aa-2bx+cxx}-a}{c}
\]
wodurch, wenn wir nach der Integration $x = \frac{b\pm\sqrt{bb-aac}}{c}$ setzen, weil ja in diese Fällen $\sqrt{aa-2bx+cxx} = 0$ und $\Pi = \Delta$ wird 
\[
\int \frac{x \mathrm{d}x}{\sqrt{aa-2bx+cxx}} = \frac{b}{c}\Delta - \frac{a}{c}
\]
sein wird.
\paragraph{§5}
Wir wollen weiter setzen $s = x\sqrt{aa-2bx+cxx}$. Es wird $\mathrm{d}s = \frac{aadx-3bxdx+2cxxdx}{\sqrt{aa-2bx+cxx}}$ sein, woher man durch Integration wieder berechnet:
\[
2cf\int\frac{xx \mathrm{d}x}{\sqrt{aa-2bx+cxx}} = 3b\int\frac{x \mathrm{d}x}{\sqrt{aa-2bx+cxx}} - aa\int\frac{\mathrm{d}x}{\sqrt{aa-2bx+cxx}} + s
\]
woraus wir sofort für den Fall $\sqrt{aa-2bx+cxx} = 0$ ableiten
\[
\int\frac{xx \mathrm{d}x}{\sqrt{aa-2bx+cxx}} = \frac{3bb-aac}{2ac}\Delta - \frac{3ab}{2cc}
\]
\paragraph{§6}
Um zu höheren Potenzen aufzusteigen, wollen wir $s = xx\sqrt{aa-2bx+cxx}$, und weil daher 
\[
ds = \frac{2aax\mathrm{d}x-5bxx\mathrm{d}x+3cx^3\mathrm{d}x}{\sqrt{aa-2bx+cxx}}
\]
wird
\[
3c\int\frac{x^3 \mathrm{d}x}{\sqrt{aa-2bx+cxx}} = 5b\int\frac{xx \mathrm{d}x}{\sqrt{aa-2bx-cxx}} - 2aa\int\frac{x \mathrm{d}x}{\sqrt{aa-2bx+cxx}} + s
\]
sein und daher weiter für den Fall, indem nach der Integration $x = \frac{b\pm\sqrt{bb-aac}}{c}$ gesetzt wird, hat man
\begin{align*}
\int \frac{x^3 dx}{\sqrt{aa-2bx+cxx}} &= \frac{5b^3-3aabc}{2c^3}\Delta - \frac{15abb}{6c^3} + \frac{2a^3}{3cc}\\
&= \left(\frac{5b^3}{2c^3} - \frac{3aab}{2cc}\right)\Delta - \frac{5abb}{2c^3} + \frac{2a^3}{3cc}
\end{align*}
\paragraph{§7}
Auf ähnliche Weise sei $s = x^3\sqrt{aa-2bx+cxx}$ und weil daher wird
\[
\mathrm{d}s = \frac{3aaxx\mathrm{d}x-7bx^3\mathrm{d}x+4cx^4\mathrm{d}x}{\sqrt{aa-2bx+cxx}}
\]
wird durch Integrieren wiederum
\[
4c\int\frac{x^4 dx}{\sqrt{aa-2bx+cxx}} = 7b\int\frac{x^3 dx}{\sqrt{aa-2bx+cxx}} \\
- 3aa\int\frac{xx dx}{\sqrt{aa-2bx+cxx}} + s
\]
dann haben wir daher für den Fall, in dem $\sqrt{aa-2bx+cxx} = 0$ wird,
\[
\int\frac{x^4 \mathrm{d}x}{\sqrt{aa-2bx+cxx}} = \left(\frac{35b^4}{8c^4} - \frac{15aabb}{4c^3} + \frac{3a^4}{8cc} \right)\Delta - \frac{35ab^3}{8c^4} + \frac{55c^3b}{24c^3}
\]
\paragraph{§8}
Damit aber die Struktur bei diesen Formel besser ausgekundschaftet werden kann, wollen wir die einzelnen Formeln durch Faktoren, wie sie durch die Strukturen entstehen, ohne jegliche Vereinfachung herbeischaffen und auf diese Weise die gefundenen Integralformen präsentieren:
\begin{flalign*}
\int\frac{\mathrm{d}x}{\sqrt{aa-2bx+cxx}} =& \Delta &\\
\int\frac{x \mathrm{d}x}{\sqrt{aa-2bx+cxx}} =& \frac{b}{c}\Delta - \frac{a}{c} &\\
\int\frac{xx \mathrm{d}x}{\sqrt{aa-2bx+cxx}} =& \left(\frac{1\cdot 3 \cdot bb}{1\cdot 2\cdot cc} - \frac{aa}{1\cdot 2\cdot c}\right)\Delta - \frac{1\cdot 3\cdot ab}{1\cdot 2\cdot cc} &\\
\int\frac{x^3 \mathrm{d}x}{\sqrt{aa-2bx+cxx}} =& \left(\frac{1\cdot 3\cdot 5b^3}{1\cdot 2\cdot 3c^3} - \frac{1\cdot 3\cdot 3aab}{1\cdot 2\cdot 3cc}\right)\Delta - \frac{1\cdot 3\cdot 5abb}{1\cdot 2\cdot c^3} + \frac{1\cdot 2\cdot 2\cdot 2a^3}{1\cdot 2\cdot 3cc} &\\
\int\frac{x^4 \mathrm{d}x}{\sqrt{aa-2bx+cxx}} =& \left( \frac{1\cdot 3\cdot 5\cdot 7b^4}{1\cdot 2\cdot 3\cdot 4c^4} - \frac{1\cdot 3\cdot 5\cdot 6aabb}{1\cdot 2\cdot 3\cdot 4c^3} + \frac{1\cdot 3\cdot 3b^4}{1\cdot 2\cdot 3\cdot 4cc}\right)\Delta &\\
& - \frac{1\cdot 3\cdot 5\cdot 7ab^4}{1\cdot 2\cdot 3\cdot 4 c^4} + \frac{1\cdot 5\cdot 11 a^3b}{1\cdot 2\cdot 3\cdot 4c^3} &
\end{flalign*}
\paragraph{§9}
Wir wollen nun im Allgemeinen die Entwicklung aufstellen, indem wir nehmen $s = x^n\sqrt{aa-2bx+cxx}$ und weil daher
\[
\mathrm{d}s = \frac{naax^{n-1}dx-(2n+1)bx^ndx + (n+1)cx^{n+1}\mathrm{d}x}{\sqrt{aa-2bx+cxx}}
\]
ist, berechnet man leicht durch Integration wiederum
\begin{flalign*}
(n+1)c\int\frac{x^{n+1}\mathrm{d}x}{\sqrt{aa-2bx+cxx}} =& (2n+1)b\int\frac{x^n \mathrm{d}x}{\sqrt{aa-2bx+cxx}}&\\
&- naa\int\frac{x^{n-1}\mathrm{d}x}{\sqrt{aa-2bx+cxx}} + x^n\sqrt{aa-2bx+cxx}
\end{flalign*}
Wenn wir aber schon vorher hervorholen: 
\[
\int\frac{x^{n-1} \mathrm{d}x}{\sqrt{aa-2bx+cxx}} = \mathcal{M}\Delta - \mathfrak{M}
\]
und 
\[
\int\frac{x^n \mathrm{d}x}{\sqrt{aa-2bx+cxx}} = \mathcal{N}\Delta - \mathfrak{N}
\]
sodass diese beiden Formeln bekannt sind, wir aus diesen so dass folgende berechnet, dass
\begin{flalign*}
\int\frac{x^{n-1} \mathrm{d}x}{\sqrt{aa-2bx+cxx}} =& \left( \frac{(2n+1)b\mathcal{N}}{(n+1)c} - \frac{naa\mathcal{M}}{(n+1)c}\right)\Delta &\\
&- \frac{(2n+1)b\mathfrak{N}}{(n+1)c} + \frac{naa\mathfrak{M}}{(n+1)c}
\end{flalign*}
Durch diese Art lassen sich die Integrationen, so weit wie es beliebt, fortsetzen, so dass all diese Integrale entweder vom Logarithmen oder von Kreisbögen abhängen, je nachdem, ob der Koeffizient entweder positiv oder negativ war. Es ist aber klar, dass Werte nicht angegeben werden können, wenn der Exponent $n$ nicht eine ganze Zahl ist.
\paragraph{§10}
Aus dieser gefundenen Integralformel, wenn nach der Integration $x = \frac{b\pm\sqrt{bb-aac}}{c}$ gesetzt wird, eoher $s=0$ wird, wird
\begin{flalign*}
naa\int\frac{x^{n-1} \mathrm{d}x}{\sqrt{aa-2bx+cxx}} =& (2n+1)b\int\frac{x^n \mathrm{d}x}{\sqrt{aa-2bx+cxx}} &\\
&- (n+1)c\int\frac{x^{n+1} \mathrm{d}x}{\sqrt{aa-2bx+cxx}}
\end{flalign*}
sein; woher, wenn wir kürzer setzen
\begin{flalign*}
\int\frac{x^{n-1} \mathrm{d}x}{\sqrt{aa-2bx+cxx}} = \mathcal{P}, \int\frac{x^n \mathrm{d}x}{\sqrt{aa-2bx+cxx}} = \mathcal{Q}, &\\
\int\frac{x^{n+1} \mathrm{d}x}{\sqrt{aa-2bx+cxx}} = \mathcal{R}, \int\frac{x^{n+2} \mathrm{d}x}{\sqrt{aa-2bx+cxx}} = \mathcal{S} &\\
etc.
\end{flalign*}
die Größen $\mathcal{P}$, $\mathcal{Q}$, $\mathcal{R}$, $\mathcal{S}$, etc. so voneinander abhängen, dass
\begin{align*}
naaP &= (2n+1)bQ - (n+1)cR \\
(n+1)aa\mathcal{Q} &= (2n+3)b\mathcal{R} - (n+2)c\mathcal{S} \\
(n+2)aa\mathcal{R} &= (2n+5)b\mathcal{S} - (n+3)c\mathcal{T} \\
(n+3)aa\mathcal{S} &= (2n+7)b\mathcal{T} - (n+4)c\mathcal{U} \\
(n+4)aa\mathcal{T} &= (2n+9)b\mathcal{U} - (n+5)c\mathcal{W} \\
etc.
\end{align*}
ist. Aus den Relationen werden folgende Rechnungen gefolgert:
\begin{align*}
\frac{\mathcal{P}}{\mathcal{Q}} &= \frac{(2n+1)b}{naa} - \frac{(n+1)c}{naa\mathcal{Q}:\mathcal{R}} \\
\frac{\mathcal{Q}}{\mathcal{R}} &= \frac{(2n+3)b}{(n+1)aa} - \frac{(n+2)c}{(n+1)aa\mathcal{R}:\mathcal{S}} \\
\frac{\mathcal{R}}{\mathcal{S}} &= \frac{(2n+5)b}{(n+2)aa} - \frac{(n+3)c}{(n+2)\mathcal{S}:\mathcal{T}} \\
\frac{\mathcal{S}}{\mathcal{T}} &= \frac{(2n+7)b}{(n+3)aa} - \frac{(n+4)c}{(n+2)\mathcal{T}:\mathcal{U}} \\
etc.
\end{align*}
daher ist klar, dass die einzelnen Brüche $\frac{\mathcal{P}}{\mathcal{Q}}$, $\frac{\mathcal{Q}}{\mathcal{R}}$, $\frac{\mathcal{R}}{\mathcal{S}}$ etc. durch die folgenden hinreichend angenehm bestimmt werden.
\paragraph{§11}
Wenn also, wie es beliebt, die gefundenen Werte dieser Ausdrücke sukzessiv ineinander eingesetzt werden, erhalten wir für den Bruch $\frac{\mathcal{P}}{\mathcal{Q}}$ einen ins Unendliche fortschreitenden Kettenbruch, der 
\[
naa\frac{\mathcal{P}}{\mathcal{Q}} = (2n+1)b - \cfrac{(n+1)^2aac}{(2n+3)b - \cfrac{(n+2)^2aac}{(2n+5)b - \cfrac{(n+3)^2aac}{(2n+7)b - \cfrac{(n+4)^2aac}{(2n+9)-etc}}}}
\]
sein wird, und so gelangen wir zu einem hinreichend gefälligem und mit offensichtlicher Struktur voranschreitenden Kettenbruch, dessen Wert immer entweder durch Logarithmen (wenn $c>0$ war) oder durch Kreisbögen (wenn $c<0$ war) ausgedrückt werden kann.
\paragraph{§12}
Wir wollen $n=1$ setzen und so wird
\[
P = \int\frac{\mathrm{d}x}{\sqrt{aa-2bx+cxx}} = \Delta
\]
sein und
\[
Q = \int\frac{x \mathrm{d}x}{\sqrt{aa-2bx+cxx}} = \frac{b}{c}\Delta - \frac{a}{c}
\]
welcher Fall uns den folgenden Kettenbruch verschafft:
\[
\frac{aac\Delta}{b\Delta - a} = 3b - \cfrac{4aac}{5b - \cfrac{9aac}{7b - \cfrac{16aac}{9b - \cfrac{25aac}{11b-etc.}}}}
\]
welcher wegen der Eleganz jeder Aufmerksamkeit würdig ist. Es hilft, hier aber zu bemerken, wenn $c$ eine negative Zahl war, dass dann alle Zähler in dem Bruch positiv werden.
\paragraph{§12[a]}
Der Kettenbruch aber scheint am Anfang wie verstümmelt; wenn daher oben ihm das $b-aac$ hinzugefügt wird, wird ein noch gefälligerer und leichterer Wert zurückgegeben. Wenn nämlich der Bruch zur Kürze mit dem Buchstaben $\mathcal{S}$ bezeichnet wird, so dass $S = \frac{aac\Delta}{b\Delta - a}$, so wird durch Hinzufügen des Gliedes dessen Wert $b - \frac{aac}{\mathcal{S}} = \frac{a}{\Delta}$ und so haben wir
\[
\frac{a}{\Delta} = b - \cfrac{aac}{3b - \cfrac{4aac}{5b - \cfrac{9aac}{7b - \cfrac{16aac}{9b - \cfrac{2baac}{11b - etc.}}}}}
\]
welcher Ausdruck bemerkenswert ist, weil bis jetzt kein Weg offen steht, auf dem der Wert des Kettenbruchs a priori gefunden werden kann.
\paragraph{§13}
Wir wollen abseits davon zwei oben bemerkte Fälle entwickeln, und die sorgfältig voneinander getrennt werden sollten. Es sei daher zuerst $c = ff$ und wir haben oben gefunden, dass 
\[
\Delta = \frac{1}{f}\log{\frac{\sqrt{(bb-aaff)}}{af-b}}
\] 
sein wird, wie dass Wurzelzeichen zweideutig angenommen werden kann. Vor allem ist es daher nötig, sodass $bb > aaff$ ist, weil andererseits der Ausdruck ins imaginäre geht; es bieten sich also zwei Fälle an, je nachdem, ob $b$ entweder eine positive oder negative Größe war. Im ersten Fall für $b>0$ und daher $b>af$ ist klar, dass der Wurzel das Zeichen $-$ gegeben werden muss, sodass
\[
\Delta = \frac{1}{f}\log{\frac{\sqrt{(bb-aaff)}}{b-af}} = \frac{1}{2f}\log{\frac{b+af}{b-af}}
\]
wird und dass wir diese Summation haben
\[
\frac{2af}{\log{\frac{b+af}{b-af}}} = b - \cfrac{aaff}{3b - \cfrac{4aaff}{5b - \cfrac{9aaff}{7b - \cfrac{16aaff}{9b - etc.}}}}
\]
woher, weil $\frac{b+af}{b-af} > 1$ ist, klar ist, dass der Wert dieses Ausdrucks positiv sein wird.
\paragraph{§14}
Wenn aber $b$ eine negative Zahl war oder wenn anstelle von $b$ $-b$ geschrieben wird, muss daher auch $b>af$ sein, dann wird $\Delta = \frac{1}{2f}\log{\frac{b-af}{b+af}}$ sein, welcher Logarithmus also negativ sein wird, oder $\Delta = -\frac{1}{2f}\log{\frac{b-af}{b+af}}$, woher man die folgende Gleichung erhält:
\[
\frac{-2af}{\log{\frac{b+af}{b-af}}} = -b-\cfrac{aaff}{-3b-\cfrac{4aaff}{-5b-\cfrac{9aaff}{-7b-\cfrac{16aaff}{-9b-etc.}}}}
\]
oder durch Verändern der Zeichen
\[
\frac{2af}{\log{\frac{b+af}{b-af}}} = b + \cfrac{aaff}{-3b + \cfrac{4aaff}{5b + \cfrac{9aaff}{-7b + \cfrac{16aaff}{9b + etc.}}}}
\]
dessen Kettenbruch also die gleiche Summe hat wie jener, die wir im vorangehenden Paragraphen gefunden haben. Diese Gleichheit aber wird durch die baldigen Rechnungen gestärkt werden.
\paragraph{§15}
Wir wollen auf selbe Art den Fall für $c=-gg$ entwickeln, für welchen wir oben gefunden haben $\Delta = \frac{1}{g}\arcsin{\frac{ag}{\sqrt{bb+aagg}}}$, welcher Wert durch den Cosinus ausgedrückt gibt $\Delta = \frac{1}{g}\arccos{\frac{b}{\sqrt{bb+aagg}}}$, woher klar ist, dass durch den Tangens der Wert noch einfacher sein wird; es wird natürlich $\Delta = \frac{1}{g}\arctan{\frac{ag}{b}}$ sein, weshalb die Summation für diesen Fall liefert
\[
\frac{ag}{\arctan{\frac{ag}{b}}} = b + \cfrac{aagg}{3b + \cfrac{4aagg}{5b + \cfrac{9aagg}{7b + \cfrac{16aagg}{9b + etc.}}}}
\]
wo keine weitere Begrenzung nötig ist.
\section*{Über Kettenbrüche, die von Logarithmen abhängen}
\paragraph{§16}
Wir wollen nun sogar weitere andere Spezialfälle, die in jeder von beiden Formen enthalten sind, anfangen, weil wir ja schon beobachtet habe, dass die Formen in §13 und §14 sich überdecken, wollen wir die erste benutzen, die war
\[
\frac{2af}{\log{\frac{b+af}{b-af}}} = b - \cfrac{aaff}{3b - \cfrac{4aaff}{5b - \cfrac{9aaff}{7b - etc.}}}
\]
und wollen den Fall betrachten, in dem $b=af$ ist, welcher jedenfalls die Summe des Bruchs gibt
\[
\frac{2af}{\log{b+af}{b-af}} = 0 = b - \cfrac{bb}{3b - \cfrac{4bb}{5b - \cfrac{9bb}{7b - etc.}}}
\]
welche durch Reduktion leicht übergeht in diesen
\[
0 = 1 - \cfrac{1}{3 - \cfrac{4}{5 - \cfrac{9}{7 - \cfrac{16}{9 - etc.}}}}
\]
\paragraph{§17}
In dieser Form dem Nichts gleichen Formel ist es notwendig, dass der Nenner des ersten Bruches $= 1$ ist und daher

\begin{minipage}[c]{0.4\textwidth}
\[
1 = 3 - \cfrac{4}{5 - \cfrac{9}{7 - etc.}}
\]
\end{minipage}
\hfill
\begin{minipage}[c]{0.4\textwidth}
\[
0 = 2 - \cfrac{9}{7 - \cfrac{16}{9 - etc.}}
\]
\end{minipage}
\vspace*{2mm}

Hier ist daher wegen derselben Methode notwendig, dass der erste Nenner $=2$ wird, sodass

\begin{minipage}[c]{0.4\textwidth}
\[
2 = 5 - \cfrac{9}{7 - \cfrac{16}{9 - etc.}}
\]
\end{minipage}
\hfill
\begin{minipage}[c]{0.4\textwidth}
\[
0 = 3 - \cfrac{9}{7 - \cfrac{16}{9 - etc.}}
\]
\end{minipage}
\vspace*{2mm}

Hier wiederum muss der erste Nenner $=3$ sein und daher

\begin{minipage}[c]{0.4\textwidth}
\[
3 = 7 - \cfrac{16}{9 - \cfrac{25}{11 - etc.}}
\]
\end{minipage}
\hfill
\begin{minipage}[c]{0.4\textwidth}
\[
0 = 4 - \cfrac{16}{9 - \cfrac{25}{11 - etc.}}
\]
\end{minipage}
\vspace*{2mm}

Wieder muss der erste Nenner $=4$ sein, sodass $4 = 9 - \frac{25}{11 - etc.}$, und auf diese Weise ist klar, dass diese Relation mit derselben Struktur im Unendlichen hat, worin das Kriterium für die Wahrheit dieser Gleichung liegt.
\paragraph{§18}
Da ja in dieser Form die Zalh $b$ größer sein muss als $af$, wollen wir nun $b = 2af$ setzen und erreichen folgende Summation
\[
\frac{2af}{\log{3}} = 2af - \cfrac{aaff}{6af - \cfrac{4aaff}{10af - \cfrac{9aaff}{14af - etc.}}}
\] 
die zurückgeführt werden zu dieser Form mit lediglich Zahlen
\[
\frac{2}{\log{3}} = 2 - \cfrac{1}{6 - \cfrac{4}{10 - \cfrac{9}{14 - \cfrac{16}{18 - etc.}}}}
\]
\paragraph{§19}
Auf ähnliche Weise können alle Buchstaben aus der Rechnung herausgeworfen werden, wenn für $b$ ein Vielfaches von $af$ genommen wird. Es sei nämlich im Allgemeinen $b = naf$ und das liefert
\[
\frac{2af}{\log{\frac{n+1}{n-1}}} = naf - \cfrac{aaff}{3naf - \cfrac{4aaff}{5naf - \cfrac{9aaff}{7naf - etc.}}}
\]
welcher Bruch zu folgender Form zurückgeführt wird
\[
\frac{2}{\log{\frac{n+1}{n-1}}} = n - \cfrac{1}{3n - \cfrac{4}{5n - \cfrac{9}{7n - etc.}}}
\]
woraus man einsieht, auf welche Weise sich alle Logarithmen durch Kettenbrüche ausdrücken lassen.
\paragraph{§20}
Es können für die Zahl $n$ gebrochene Zahlen genommen werden, dann liefern aber die ersten Terme bei den einzelnen Gliedern gebrochenes, welche sich durch Reduktion zu ganzen zurückführen lassen; aber Fälle dieser Art können leicht aus der allgemeinen Form berechnet werden, in dem man sofort $b=n$ und $af=m$ schreibt; dann haben wir nämlich
\[
\frac{2m}{\log{\frac{n+m}{n-m}}} = n - \cfrac{mm}{3n - \cfrac{4mm}{5n - \cfrac{9mm}{7n - etc.}}}
\]
woher, wenn anstelle von $m$ $\sqrt{k}$ geschrieben wird
\[
\frac{2\sqrt{k}}{\log{\frac{n+\sqrt{k}}{n-\sqrt{k}}}} = n - \cfrac{k}{3n - \cfrac{4k}{5n - \cfrac{9k}{7n - etc.}}}
\]
sein wird.
\paragraph{§21}
Daher können die hyperbolischen Logarithmen der ganzen Zahlen durch Kettenbrüche ausgedrückt werden. Es sei also im allgemeinen die ganze Zahl $i$ vorgelegt und man setze $\frac{n+m}{n-m} = i$, es wird $\frac{n}{m} = \frac{i+1}{i-1}$ sein. Man nimmt also $n=i+1$ und $m=i-1$ und wir haben
\[
\frac{2(i-1)}{\log{i}} = i + 1 - \cfrac{(i-1)^2}{3(i+1) - \cfrac{4(i-1)^2}{5(i+1) - \cfrac{9(i-1)^2}{7(i+1) - \cfrac{16(i-1)^2}{9(i+1) - etc.}}}}
\]
woher wir berechnen
\[
\log{i} = \cfrac{2(i-1)}{i+1 - \cfrac{(i-1)^2}{3(i+1) - \cfrac{4(i-1)^2}{5(i-1) - \cfrac{9(i-1)^2}{7(i+1) - etc.}}}}
\]
\paragraph{§22}
Wenn wir Brüche solcher Art für Logarithmen gebrochener Zahlen haben wollen, sollten wir $\frac{n+m}{n-m} = \frac{p}{q}$ setzen, woher $n = p + q$ und $m = p - q$ ist, weshalb wir haben
\[
\log{\frac{p}{q}} = \cfrac{2(p-q)}{1(p+q) - \cfrac{1(p-q)^2}{3(p+q) - \cfrac{4(p-q)^2}{5(p+q) - \cfrac{9(p-q)^2}{7(p+q) - etc.}}}}
\]
welcher Form umso bemerkenswerter ist, weil hinreichend angehmer angewendet können, um Näherungen für Logarithmen zu finden. Je mehr aber der Kettenbruch konvergiert, desto kleiner war der Bruch $\frac{p-q}{p+q}$.
\paragraph{§23}
Um ein Beispiel zu illustrieren, wollen wir $p=2$ und $q=1$, woher freilich nicht allzu schnelle Konvergenz erwartet werden darf, und es wird
\[
\log{2} = \cfrac{2}{3 - \cfrac{1}{9 - \cfrac{4}{15 - \cfrac{9}{21 - etc.}}}}
\]
sein, indem nur das erste Glied $2/3$ in Dezimalstellen genommen wird, es $0.666666$ liefert, während aus Tabellen der Wert $\log{2} = 0.693147$ erhalten wird, wo der Fehler schon klein genug ist. Wollen wir schon die beiden ersten Glieder nehmen $\frac{2}{3 - \frac{1}{9}} = \frac{9}{13} = 0.6923$. Durch Nehmen dreier Glieder aber haben wir 
\[
\cfrac{2}{3- \cfrac{1}{9 - \cfrac{4}{15}}} = \cfrac{2}{3 - \cfrac{15}{131}} = \cfrac{262}{378} = 0.693121
\]
welcher Wert von der Wahrheit nur um die $0.000026$ verschieden ist. Viel schneller zeigt sich aber die Konvergenz, wenn wir $p=3$ nehmen und $q=2$, sodass wir haben
\[
\log{\cfrac{3}{2}} = \cfrac{2}{5 - \cfrac{1}{15 - \cfrac{4}{25 - \cfrac{9}{35 - etc.}}}}
\]
dessen erstes Glied $\frac{2}{5} = 0.400000$ gibt; in Wirklichkeit ist aber $\log{\frac{3}{2}} = 0.405465108$. Durch $2$ genommene Glieder $\frac{2}{5- \frac{1}{15}}$ wird $\log{\frac{3}{2}} = 0.40540$ berechnet, wo der Fehler sich erst in der fünften Stelle einschleicht. Es werden drei Glieder genommen
\[
\cfrac{2}{5 - \cfrac{1}{15 - \cfrac{4}{25}}} = \cfrac{2}{5 - \cfrac{25}{371}} = 0.4054654
\]
wo der Fehler sich schließlich in der siebten Stelle bemerkbar macht.
\paragraph{§24}
Wegen dieses ausgezeichneten Nutzens, der sich entgegen der Erwartung gezeigt hat, wird es der Mühe wert sein, solch Untersuchungen ins Allgemeine auszuweiten; und für dieses Ziel wollen wir die Formel zwischen den Buchstaben $m$ und $n$ -- oben in §20 gegeben -- benutzen, wo
\[
\log{\frac{n+m}{n-m}} = \cfrac{2m}{n - \cfrac{m}{3n - \cfrac{4mm}{5n - \cfrac{9mm}{7n - \cfrac{16mm}{9n - etc.}}}}}
\]
ist, woher, wenn wir nur das erste Glied nehmen, fast $\log{\frac{n+m}{n-m}} = \frac{2m}{n}$ sein wird, durch die ersten beiden Glieder genommen
\[
\cfrac{2m}{n - \cfrac{3mm}{3n}}
\]
wird schon näher $\log{\frac{n+m}{n-m}} = \frac{6mn}{3nn-mm}$ sein; durch die ersten drei Glieder genommen wird
\[
\log{\frac{n+m}{n-m}} = \cfrac{2m}{n- \cfrac{mm}{3n - \cfrac{4mm}{5n}}} = \cfrac{30mnn - 8m^3}{15n^3 - 9mmn}
\]
sein.
\paragraph{§25}
Es ist nicht so sehr nötig, diese Brüche weiter fortzusetzen, wir wollen den schon gefundenen Brüchen nämlich den Bruch $\frac{0}{1}$ vorne anheften, dass wir diese Progression der Brüche erhalten
\begin{center}
\begin{tabular}[c]{cccc}
$\mathrm{I}$ & $\mathrm{II}$ & $\mathrm{III}$ & $\mathrm{IV}$  \\[2mm]
$\dfrac{0}{1}$ & $\dfrac{2m}{n}$ & $\dfrac{6nn}{3nn-mm}$ & $\dfrac{30mnn-8m^3}{15n^3-9mmn}$
\end{tabular}\\[4mm]
\end{center}
dessen Zähler wie Nenner aus zwei vorangehenden zur Einfachheit einer rekurrenten Reihe geformt werden können. Der dritte ist natürlich aus dem ersten und zweiten mit Hilfe der Beziehungsskala $3n, -mm$ geformt werden können; der vierte aber aus den beiden vorangehenden mit Hilfe der Beziehungsskala $5n, -4mm$. Für den fünften aber ist die Beziehungsskala $7n, -9mm$ zu benutzen, für den sechsten $9n, -16mm$ und so weiter. Daher findet man auf diese Weise leicht den fünften Bruch
\begin{center}
\begin{tabular}[c]{c}
$\mathrm{V}$ \\[2mm]
$= \dfrac{210mn^3 - 111m^3n}{105n^4 - 90mmnn + 9m^4}$
\end{tabular}\\[4mm]
\end{center}
auf ähnliche Weise
\begin{center}
\begin{tabular}[c]{c}
$\mathrm{VI}$\\[2mm]
$\dfrac{1980n^4 - 1470m^3nn + 128m^5}{945n^5 - 1050mmn^3 + 225m^4n}$\\[4mm]
etc.
\end{tabular}\\[4mm]
\end{center}
\paragraph{§26}
Es hilft aber ganz besonders zu bemerken, dass diese Brüche stetig länger werden und sich durch Zunahmen weniger an die Wahrheit annähern. Die Zuwächse schreiten mit der Struktur voran, wie hier zu sehen ist:
\begin{flalign*}
\mathrm{II} - \mathrm{I} &= \frac{2m}{n} &\\
\mathrm{III} - \mathrm{II} &= \frac{2m^3}{n(3nn-mm)} &\\
\mathrm{IV} - \mathrm{III} &= \frac{2\cdot 4m^5}{(3nn-mm)(15n^3-9mmn)} &\\
\mathrm{V} - \mathrm{IV} &= \frac{2\cdot 4\cdot 9m^4}{(15n^3-9mmn)(105n^4 - 9mmnn+9m^4)} &\\
\mathrm{VI} - \mathrm{V} &= \frac{2\cdot 4\cdot 9\cdot 16m^9}{(105n^4 - 90mmnn + 9m^4)(945n^5-1050mmn^3+225m^4n)} &
\end{flalign*}
woher klar ist, je größer die Zahl $n$ in Beuzg auf $m$ war, dass die Differenzen umso schneller klein werden, dass sie ohne Fehler missachtet werden können.
\section*{Über Kettenbrüche, die von Kreisbögen abhängen}
\paragraph{§27}
Der Kreisbogen aus §15, dessen Tangens $\frac{ag}{b}$ ist, wird durch einen Bruch so ausgedrückt, dass
\[
\arctan{\frac{ag}{b}} = \cfrac{ag}{b + \cfrac{aagg}{3b + \cfrac{4aagg}{5b + \cfrac{9aagg}{7b +etc.}}}}
\]
ist. Wir wollen nun zur Ähnlichkeit zur oberen Formen $ag=m$ und $b=n$ setzen und haben
\[
\arctan{m}{n} = \cfrac{m}{n + \cfrac{mm}{3n + \cfrac{4mm}{5n + \cfrac{9mm}{7n + etc.}}}}
\]
welche Formel umso schneller konvergiert, je größer $n$ in Bezug auf $m$ war; daher ist klar, dass auch dieser Ausdruck ertragreich bei der Rechnung angewendet werden kann.
\paragraph{§28}
Wir wollen mit dem Fall, in dem $m=1$ und $n=1$ ist, anfangen, wodurch
\[
\arctan{\frac{m}{n}} = \frac{\pi}{4} = \cfrac{1}{1 + \cfrac{1}{3 + \cfrac{4}{5 + \cfrac{9}{7 + \cfrac{16}{9 + etc.}}}}}
\]
ist, welcher Bruch freilich nicht so sehr konvergiert; Aber dennoch sehen wir, wie er allmählich sich der Wahrheit nähert, da wir ja wissen, dass $\frac{\pi}{4} = 0.78539816339$ ist.
Und das erste Glied gibt
\[
\frac{\pi}{4} = \frac{1}{1} \text{(zu groß)}
\]
zwei Glieder liefern
\[
\frac{\pi}{4} = \cfrac{1}{1 + \cfrac{1}{3}} = \frac{3}{4} \text{(zu klein)}
\]
drei Glieder liefern
\[
\frac{\pi}{4} = \cfrac{1}{1 + \cfrac{1}{3 + \cfrac{4}{5}}} = \frac{19}{24} = 0.7916 \text{(0.7816)}
\]
Man sollte vier Glieder nehmen, sodass
\[
\frac{\pi}{4} = \cfrac{1}{1 + \cfrac{1}{3 + \cfrac{4}{5 + \cfrac{9}{7}}}}
\]
wird, wo der Fehler schließlich in der dritten Stelle gefunden wird.
Im Übrigen ist dieser Kettenbruch diesem so sehr ähnlich, den einst Brouncker vorgebracht hat, welcher sich so verhält:
\[
\cfrac{1}{2 + \cfrac{9}{2 + \cfrac{25}{2 + \cfrac{49}{2 + etc.}}}}
\]
\paragraph{§29}
Damit wir aber einen konvergenteren Kettenbruch erreichen, wollen wir $\arctan{\frac{m}{n}} = 30$° setzen, weil dessen Tangens $\frac{1}{\sqrt{3}}$ ist; damit $n$ nicht eine irrationale Zahl ist, wollen wir $m=\sqrt{3}$ und $n=3$ setzen; daher wird also
\[
\frac{\pi}{6} = \cfrac{\sqrt{3}}{3 + \cfrac{3}{9 + \cfrac{12}{15 + \cfrac{27}{21 + \cfrac{48}{27 + etc.}}}}}
\]
welcher Form zurückgeführt wird zum folgenden
\[
\frac{\pi}{6\sqrt{3}} = \cfrac{1}{3 + \cfrac{1}{3 + \cfrac{4}{15 + \cfrac{9}{7 + \cfrac{16}{27 + \cfrac{25}{11 + etc.}}}}}}
\]
für welche Entwicklung wir zuerst den Näherungswert $\frac{\pi}{6\sqrt{3}}$ suchen wollen, der $0.3022998$ ist. Nun liefert aber das erste Glied $\frac{\pi}{6\sqrt{3}} = 0.3333$, zwei Glieder liefern
\[
\frac{\pi}{6\sqrt{3}} = \cfrac{1}{3 + \cfrac{1}{3}} = \frac{3}{10} = 0.3000
\]
und drei Glieder liefern
\[
\frac{\pi}{6\sqrt{3}} = \cfrac{1}{3 + \cfrac{1}{3 + \cfrac{4}{15}}} = \frac{49}{102} = 0.30247
\]
wo der Fehler die vierte Stelle betrifft.
\paragraph{§30}
Es kann aber viel schneller Konvergenz verschafft werden, wenn wir den rechten Winkel in zwei Teile teilen, weil ich ja einst gesetzt habe, dass $\arctan{\frac{1}{2}}+\arctan{\frac{1}{3}} = \arctan{1} = \frac{\pi}{4}$ ist.
Und so finden wir zwei Kettenbrüche, deren Summe den Wert $\frac{\pi}{4}$ gibt, die da sind

\begin{minipage}[c]{0.5\textwidth}
\[
\arctan{\frac{1}{2}} = \cfrac{1}{2 + \cfrac{1}{6 + \cfrac{4}{10 + \cfrac{9}{14 + etc.}}}}
\]
\end{minipage}
\begin{minipage}[c]{0.5\textwidth}
\[
\arctan{\frac{1}{3}} = \cfrac{1}{3 + \cfrac{1}{9 + \cfrac{4}{15 + \cfrac{9}{21 + etc.}}}}
\]
\end{minipage}
\vspace*{2mm}

Es ist aber klar, dass beide Brüche und besonders der hintere sehr stark konvergieren.
\paragraph{§31}
Wir wollen nun aber unseren allgemeinen Kettenbruch in einen gewöhnlichen Bruch umwandeln; und aus dem ersten Glied allein finden wir $\arctan{\frac{m}{n}} = \frac{m}{n}$, zwei Glieder liefern $\arctan{\frac{m}{n}} = \frac{3mm}{3nn+mm}$, drei Glieder liefern $\arctan{\frac{m}{n}} = \frac{15mnn+4m^3}{15n^3+9mmn}$.
Es werden vier Glieder genommen, woraus $\arctan{\frac{m}{n}} = \frac{105mn^3 + 55m^3n}{105n^4 + 90mmnn+9m^4}$ wird.
Wenn daher diesen Brüchen nun wie oben $\frac{0}{1}$ vorangestellt wird, entsteht diese Progression
\begin{center}
\begin{tabular}[c]{ccccc}
$\mathrm{I}$ & $\mathrm{II}$ & $\mathrm{III}$ & $\mathrm{IV}$ & $\mathrm{V}$ \\[1mm]
$\dfrac{0}{1}$ & $\dfrac{m}{n}$ & $\dfrac{3mn}{3nn+mm}$ & $\dfrac{15mnn+4m^3}{15n^3+9mmn}$ & $\dfrac{105mn^3 + 55m^3n}{105n^4 + 90mmnn+9m^4}$
\end{tabular}\\[4mm]
\end{center}
dessen einzelne Terme genauso aus den zwei vorangegangenen nach einem festen Gesetz gebildet werden können; dann ist natürlich
\begin{flalign*}
& \text{für } \mathrm{III} \text{ die Beziehungsskala } 3n, +mm &\\
& \text{für } \mathrm{IV}  \text{ die Beziehungsskala } 5n, +4mm &\\
& \text{für } \mathrm{V}   \text{ die Beziehungsskala } 7n, +9mm &\\
\end{flalign*}
\end{document}